\documentclass{article}
\usepackage{mathptmx}      % use Times fonts if available on your TeX system

\usepackage[utf8]{inputenc}
\usepackage[british]{babel}
\usepackage[T1]{fontenc}
\usepackage{amsfonts}
\usepackage{color}
\usepackage{bm}
\usepackage{amsmath,amsthm}

\newtheorem{theorem}{Theorem}[section]
\newtheorem{lemma}[theorem]{Lemma}
\newtheorem{proposition}[theorem]{Proposition}

\newtheorem{example}[theorem]{Example}

\newcommand{\half}{\frac{1}{2}}

\newcommand{\inv}[1]{\frac{1}{#1}}

\def\R{\mathbb{R}}
\def\N{\mathbb{N}}

\def\S{\mathbb{S}}
\def\0{{\bf 0}}
\def\x{{\bf x}}
\def\y{{\bf y}}
\def\rE{\mathrm{E}}     %expectation,

\def\S{\mathbb{S}}
\def\D{\mathbb{D}}
\def\0{{\bf 0}}
\def\x{{\bf x}}
\def\y{{\bf y}}
\def\rE{\mathrm{E}}     %expectation,

\def\Md{\mathbb{M}^d}
\def\P{\mathbb{P}}

\def\0{{\bf 0}}
\def\rE{\mathrm{E}}     %expectation,
\def\var{\mathop{\rm var}\nolimits}    %variance
\def\cov{\mathop{\rm cov}\nolimits}    %covariance

\newcommand\be{\begin{equation}}
\newcommand\ee{\end{equation}}

\newcommand{\eq}[1]{Eq. \eqref{#1}}

\newcommand\bes{\begin{eqnarray}}
\newcommand\ees{\end{eqnarray}}
\newcommand\non{\nonumber}

\newcommand{\hf}[1]{\frac{#1}{2}}

\newcommand\xx{\mathbf{x}}

\newcommand\comment[1]{{}}
\newcommand\eps{\epsilon}

\begin{document}

\title{Strong Local  Nondeterminism and Exact  Modulus of  Continuity for Isotropic Gaussian Random Fields on 
Compact Two-Point Homogeneous Spaces}

\author{Tianshi  Lu$^1$, Chunsheng Ma$^1$, Yimin Xiao$^2$\\
$^1$Dept. of Mathematics and Statistics, Wichita State University, Wichita, KS 67260, USA\\
$^2$Dept. of Statistics and Probability, Michigan State University, East Lansing, MI 48824, USA}

\maketitle

\begin{abstract}
This  paper is concerned with sample path properties of isotropic  Gaussian fields on compact two-point homogeneous spaces. 
In particular, we establish the property of strong local nondeterminism of an isotropic Gaussian field based on the high-frequency 
behaviour of its angular power spectrum, and then exploit this result to establish an exact uniform modulus of continuity for its 
sample paths. 
 
\end{abstract}

\section{Introduction}

The analysis of sample path properties of random fields has been considered by many authors \cite{FalconerXiao2014}, 
\cite{NanaXiao2012}, \cite{Pitt1978},  \cite{Xiao2007}, \cite{Xiao2013}, but the index set of the random fields is typically 
restricted to be the Euclidean space $\R^d$. Recently the investigation of  sample path properties of random fields over 
the unit sphere $\S^d$ has been conducted by  \cite{Lan2018}, \cite{LanXiao2018}, \cite{LangSchwab2015}.
This paper is concerned with sample path properties of isotropic Gaussian random fields on 
a  $d$-dimensional compact two-point homogeneous space  $\mathbb{M}^d$. 

It is well-known that $\mathbb{M}^d$ is a compact Riemannian symmetric space of rank one, and  belongs to one of the
following five families  (\cite{Helgason2011},   \cite{Wang1952}): the unit spheres $\S^d$ ($ d =1, 2, \ldots$), the real 
projective spaces $\mathbb{P}^d(\R)$ ($d = 2, 4, \ldots$), the complex projective spaces  $\mathbb{P}^d(\mathbb{C})$ 
($d = 4, 6, \ldots$), the quaternionic  projective spaces $\mathbb{P}^d(\mathbb{H})$ ($d = 8, 12, \ldots$), and the Cayley  
elliptic plane $\mathbb{P}^{16}  (Cay)$ or $\mathbb{P}^{16}  (\mathbb{O})$.
There are at least two different approaches to the subject of compact two-point homogeneous spaces   \cite{MaMalyarenko2018},
including an approach based on Lie algebras and a geometric approach, which are used in   probabilistic  and statistical literature 
\cite{Askey1976}, \cite{Gangolli1967},  \cite{Malyarenko2013},   \cite{Patrangenaru2016}.
All compact two-point homogeneous spaces share the same property  that all   geodesics in a given one of these spaces are closed 
and have the same length \cite{Gangolli1967}. In particular, when the unit sphere $\mathbb{S}^d$ is embedded into the space 
$\mathbb{R}^{d+1}$, the length of any geodesic line is equal to that of the unit circle, that is, $2\pi$. In what follows, the distance 
$\rho (\x_1, \x_2)$ between two points $\x_1$ and $\x_2$ on $\Md$ is defined in such a way that  the length of any geodesic 
line on all $\mathbb{M}^d$ is equal to $2\pi$, or the distance between any two points is bounded between 0 and $\pi$, {\em i.e.},
$0 \le \rho (\x_1, \x_2) \le  \pi$.  Over $\S^d$, for instance, $\rho (\x_1, \x_2) $ is defined by $\rho (\x_1, \x_2)= \arccos (\x_1' \x_2)$ 
for all $ \x_1, \x_2 \in \S^d$, where $\x_1'\x_2$ is the inner product between $\x_1$ and $\x_2$.
Expressions of
$\rho (\x_1, \x_2)$ on other spaces may be found in \cite{Bhattacharya2012}.

Gaussian random fields on $\Md$ have been studied  in \cite{Askey1976}, \cite{CGLP20}, \cite{Gangolli1967}, \cite{LuMa2020},
\cite{MaMalyarenko2018},  \cite{Malyarenko2013}, 
among others, while theoretical investigations and
practical applications of  scalar and vector random fields on spheres may be found in
\cite{Askey1976}, \cite{Bingham1973},    \cite{Cohen2012}, 
              \cite{Gangolli1967},
             \cite{Leonenko2012},    \cite{Ma2015}, \cite{Ma2017},
         \cite{Malyarenko2013},  \cite{Malyarenko1992},
    \cite{Yadrenko1983}-\cite{Yaglom1987}.
Recently, a series representation for a real-valued isotropic Gaussian random field  on $\Md$  is presented in \cite[Chapter 2]{Malyarenko2013}. More generally,  a series representation  is provided  in    \cite{MaMalyarenko2018} for a vector random field that is isotropic and mean square continuous on $\Md$ and stationary on a temporal domain, and a general form of the covariance matrix function is derived  for such a vector random field, which involve Jacobi polynomials and the distance defined on $\Md$.  Parametric and semiparametric covariance matrix structures on $\Md$ are constructed   in \cite{LuMa2020}.

\begin{table}
 \centering
\caption{ Parameters $\alpha$ and $\beta$ associated with Jacobi polynomials over $\Md$}
  \label{table1}
    \begin{tabular}{|l| c |  c|}
    \hline
    $\mathbb{M}^d$ &    $\alpha$ &  $\beta$ \\
    \hline
    $\mathbb{S}^d$, $d=1$, $2$, \dots &  $\frac{d-2}{2}$ &  $\frac{d-2}{2}$  \\
    $\mathbb{P}^d(\mathbb{R})$, $d=2$, $3$, \dots &   $\frac{d-2}{2}$  &   $-\frac{1}{2}$  \\
    $\mathbb{P}^d(\mathbb{C})$, $d=4$, $6$, \dots &  $\frac{d-2}{2}$  & 0  \\
    $\mathbb{P}^d(\mathbb{H})$, $d=8$, $12$, \dots &   $\frac{d-2}{2}$  &  1 \\
    $\mathbb{P}^{16}(Cay)$ &    7   &  3 \\
    \hline
  \end{tabular}
\end{table}

A second-order random field $Z = \{ Z (\x), \x \in \Md \}$ is called stationary (homogeneous) and isotropic, 
if its mean function $\rE Z (\x) $ does not depend on $\x$, and its covariance  function,
 $$ \cov (  Z (\x_1),  Z ( \x_2) ) = \rE [ ( Z (x_1) - \rE Z (\x_1)) ( Z (x_2) - \rE Z (\x_2)) ],    ~~~~~~ \x_1, \x_2 \in \Md, $$
depends only on the distance $\rho (\x_1, \x_2)$ between $\x_1$ and $\x_2$.
We denote such a covariance function by $C( \rho (\x_1, \x_2)), \, \x_1, \x_2 \in \Md, $
and call it an isotropic covariance  function on $\Md$.
An isotropic  random field $\{ Z (\x), \x \in \Md \}$ is said to be mean square continuous if
 $$ \rE  | Z (\x_1) -Z (\x_2) |^2 \to 0,  ~~ \mbox{as} ~~ \rho (\x_1, \x_2 ) \to 0,  ~ \x_1, \x_2 \in \mathbb{M}^d. $$
It implies the continuity of  of the  associated covariance  function  in terms of $\rho (\x_1, \x_2)$.

For an  isotropic and mean square continuous random field on $\Md$, its covariance  function is of the form  
(\cite{LuMa2020}, \cite{MaMalyarenko2018}) 
 \begin{equation}
 \label{cov.mf1}
 C( \rho (\x_1, \x_2))  = \sum_{n=0}^\infty  b_n   
                             \frac{P_n^{(\alpha, \beta) } \left( \cos  \rho (\x_1, \x_2) \right)}{P_n^{(\alpha, \beta) } (1)},
      ~~~ \x_1, \x_2 \in \mathbb{M}^d,
      \end{equation}
      where  $\{ b_n,  n \in \mathbb{N}_0 \}$ is a  summable sequence of nonnegative constants, 
       \begin{equation}
        \label{JacobiPolynomial}
          P_n^{(\alpha, \beta)} (x) = \frac{\Gamma (\alpha+n+1)}{n! \Gamma (\alpha+\beta+n+1)}\sum_{k=0}^n\binom{n}{k}
          \frac{\Gamma (\alpha+\beta+n+k+1)}{\Gamma ( \alpha+k+1 )} \left(\frac{x-1}{2} \right)^k,
          \end{equation}
         \hfill $ x \in \R, \quad n \in \mathbb{N}_0, $

\noindent
are the Jacobi polynomials \cite{Szego1975} with specific pair of parameters $\alpha$ and $\beta$ given in Table \ref{table1}, 
\begin{equation}\label{Eq:Fact2}
P_n^{(\alpha, \beta)}(1) = {n+\alpha \choose n} = \frac{\Gamma(\alpha+n+1)}{n! \Gamma(\alpha + 1)}, \qquad \quad ~~~~~~~~ n \in \N_0,
\end{equation}
and $\mathbb{N}_0$ and $\N$ denote the sets of nonnegative integers and of positive integers, respectively.
On the other hand, if $C( \rho (\x_1, \x_2))$ is  a function  of the form (\ref{cov.mf1}), then there exists an  
isotropic Gaussian or elliptically contoured random field on $\Md$ with $ C( \rho (\x_1, \x_2))$  as its covariance  
function  \cite{LuMa2020}, \cite{MaMalyarenko2018}.

For a    centered  isotropic  Gaussian random field  $ \{ Z(\x), \x \in \Md \}$   with covariance function $C( \rho (\x_1, \x_2))$
given by (\ref{cov.mf1}), its variogram is 
\begin{equation}\label{Eq:viog1}
  \gamma (\x_1, \x_2) = \frac{1}{2} \rE \big[(Z(\x_1)- Z(\x_2))^2\big] = \sum_{n=0}^\infty b_n \bigg(1- \frac{P_n^{(\alpha, \beta) }
   \left( \cos  \rho (\x_1, \x_2) \right)}  {P_n^{(\alpha, \beta) } (1)} \bigg)
\end{equation}
for all $ \x_1, \x_2 \in \Md$. We will see that
many of the probabilistic and regularty properties of $Z$ are determined by the asymptotic property of $\{b_n, n \ge 0\}$.

For an isotropic and mean square continuous Gaussian random field $ \{ Z(\x), $ $\x \in \Md \}$, 
Section 2 establishes the property of strong local nondeterminism in terms of an asymptotic condition on 
the angular power spectrum $\{b_n, n \in \N_0 \}$.  Section 3 determines its exact uniform modulus of continuity. 
The proofs of propositions and theorems are given in Section 4.  

 \section{Strong local nondeterminism}
 
In what follows let $d \ge 2$,  $\alpha = \frac{d-2}{2}$, and let $\beta$ be given in the last column of Table \ref{table1} 
associated with $\alpha$.  Our focus is on a Gaussian random field  $\{ Z(\x), \x \in \Md \}$ that is isotropic and mean 
square continuous on $\Md$, whose covariance function is known   (\cite{LuMa2020}, \cite{MaMalyarenko2018}) to be 
of the form  (\ref{cov.mf1}). This section establishes the property of strong local nondeterminism (SLND) for  
 $\{ Z(\x), \x \in \Md \}$ under certain asymptotic condition on the coefficient sequence  $\{ b_l, l \in \N_0 \}$ in (\ref{cov.mf1}).

Denote by $\var ( Z(\x) |  Z(\x_1), \ldots, Z(\x_n))$ the conditional variance of $Z(\x)$ given $Z(\x_1),$ $\ldots, Z(\x_n)$.
For a Gaussian random field $\{ Z(\x), \x \in \D \}$,  it is known  that
   \begin{equation}
   \label{inf.eq}
    \var (Z(\x) | Z(\x_1), \ldots, Z(\x_n)) = \inf  \,  \rE \bigg( Z(\x) - \sum_{k=1}^n a_k Z(\x_k) \bigg)^2, 
    \end{equation}
  where the infimum is taken over all $(a_1, \ldots, a_n)' \in \R^n$.

For an isotropic and mean square continuous Gaussian random field on $\Md$, the SLND property is described 
in the following theorem, under the condition that the coefficients in the Jacobi  expansion of its covariance 
function fulfill  inequality \eqref{thm1.ineq1} below.  In the particular case of $d=2$  and $\Md=\S^d$, the SLND property
was  derived in \cite{Lan2018}.

\begin{theorem}
\label{thm1}
Suppose that  $\{ Z(\x), \x \in \Md \}$ is an isotropic and mean square continuous Gaussian random field with mean 0 
and covariance function (\ref{cov.mf1}), 
where $\{ b_l, l \in \N_0 \}$ is a summable sequence with nonnegative terms.
If there are $l_0 \in \N$ and  positive constants $\nu$ and $\gamma_1$,   such that
 \begin{equation}
  \label{thm1.ineq1}
   b_l  (1+l)^{1+\nu}   \ge \gamma_1,  ~~~~~ \forall ~  l \ge l_0, 
   \end{equation} 
 then there is a positive constant  $\gamma$  such that the inequality
 \begin{equation}
 \label{thm1.ineq2}
  \var ( Z(\x) |  Z(\x_1), \ldots, Z(\x_n)) \ge  \gamma  \left(  \min\limits_{1 \le  k \le n} \rho (\x, \x_k)  \right)^\nu
 \end{equation}
holds for all  $n \in \N$ and all $ \x, \, \x_k \in \Md$  ($k=1, \ldots, n$).
  \end{theorem}
  
  Inequality  \eqref{thm1.ineq1} implies that $\{ b_l, l \in \N_0 \}$ is  away from zero for all large  $l \in \N$. 
  In this case, the covariance function (\ref{cov.mf1}) of   $\{ Z(\x), \x \in \Md \}$  is strictly positive definite 
  \cite{Barbosa2016}.

To prove Theorem \ref{thm1}, we will make use of  Propositions \ref{prop2} and \ref{prop1} below. 
Proposition  \ref{prop2} is quite interesting in its own right, since  $\D$ is an arbitrary index set.  

 \begin{proposition}
 \label{prop2}
   If $C(\x_1, \x_2)$ is a covariance function on $\D$, then so is the function $ C(\x_0, \x_0) C(\x_1, \x_2)-
   C(\x_1, \x_0) C(\x_2, \x_0)$, $\x_1, \x_2 \in \D$, for any fixed point $\x_0 \in \D$.  
   \end{proposition}
 
For a continuous  function $g(\vartheta)$ on $[0, \pi]$, we expand it in terms of the  Jacobi polynomials  
\cite{Szego1975}   as follows
 \begin{equation}     \label{fexp}
g( \vartheta ) 
= \sum_{n=0}^\infty b^{ (\alpha,\beta)}_n( g )\frac{P^{(\alpha,\beta)}_n(\cos \vartheta )}{P^{ (\alpha,\beta)}_n(1)},  
  ~~~~~ \vartheta \in [0, \pi],
 \end{equation}
where $  P^{ (\alpha,\beta ) }_n(1)$ is given by (\ref{Eq:Fact2}) and the coefficients $b^{ (\alpha,\beta)}_n (g) $  are given by 
\begin{equation}
 \label{cab}
 \begin{split}
b^{ (\alpha,\beta)}_n (g) &=\frac{(2n+\alpha+\beta+1)\Gamma(n+\alpha+\beta+1)}{\Gamma(n+\beta+1)\Gamma(\alpha+1)}\\
&\qquad \   \times \int_0^\pi g(x)P_n^{(\alpha,\beta)}(\cos x)
\sin^{2\alpha+1} \left( \frac{x}{2} \right) \cos^{2\beta+1}   \left( \frac{x}{2} \right) dx,  ~~~~  n \in \N_0.
\end{split}
\end{equation} 
 
In the proof of  Proposition \ref{prop1} we actually construct a specific function $g_\eps (\vartheta)$ that 
satisfies the following properties (i) - (iii).  Such a function is termed  as a spherical bump function \cite{Lan2018}  
 in the case of $\Md=\S^d$. 
 
 \begin{proposition}
 \label{prop1}
 For three constants $r > 1$, $n_0 \in \N$,  and $ \eps \in (0, \pi]$,  there exist  a continuous function 
 $g_\eps (\vartheta)$ on $[0, \pi]$ and a positive constant $\gamma_r$, such that
\begin{enumerate}
\item[\textup{(i)}]  $g_\eps(0)=1$, and $g_\eps(\pi)=0$,
\item[\textup{(ii)}]  $b^{(\alpha,\beta)}_n( g_\eps)=0$ \, for  $0\le n<n_0$, 
\item[\textup{(iii)}]  $|b^{(\alpha,\beta)}_n(g_\eps)|\le  M_r\eps(1+n\eps)^{-r}$.
\end{enumerate}
 \end{proposition}

 \section{Modulus of continuity}
 
 As an application of Theorem \ref{thm1}, we determine the exact uniform modulus of continuity of $Z$ 
 on $\Md$. Theorem \ref{thm2}
 improves Corollary 5.3 in \cite{CGLP20} significantly.
 
 \begin{theorem}
\label{thm2}
Assume that  $\{ Z(\x), \x \in \Md \}$ is an isotropic and mean square continuous  Gaussian random field 
with covariance function  (\ref{cov.mf1}).  
  If  there are   $l_0 \in \N$ and positive  constants $\gamma_1, \gamma_2$,  $\nu \in (0, 2)$,  such that
  \begin{equation}
  \label{thm2.ineq2}
  \gamma_1 \le   b_l  (1+l)^{1+\nu} \le \gamma_2,   ~~~~~ \forall \ l \ge l_0, 
   \end{equation} 
then there is a positive and finite constant $\kappa$  such that, with probability 1, 
 \begin{equation}
 \label{thm2.eq}
  \lim_{\varepsilon \to 0} \sup\limits_{ \substack{ \x_1, \x_2 \in \Md  \\   \rho (\x_1, \x_2) \le \varepsilon }  }
    \frac{ | Z(\x_1) -Z(\x_2) |}{  \rho (\x_1, \x_2)^{\nu/2} \sqrt{ \ln  \rho (\x_1, \x_2)} } = \kappa.
 \end{equation}
 \end{theorem}

Under the assumption   (\ref{thm2.ineq2})  on the coefficients $\{b_l,\, l \in \N_0\}$  in   (\ref{cov.mf1}), we obtain 
upper and lower bounds for the variogram of the Gaussian random field $\{ Z(\x), \x \in \Md \}$ in terms of the 
distance function over 
$\Md$ in Proposition 3, which will be employed to prove Theorem \ref{thm2}.  Some related asymptotic relationships 
between $\{b_l,  l \in \N_0  \}$ and  the variogram may be found in  \cite{Malyarenko05}.  Under a  stronger condition 
than (\ref{thm2.ineq2}), the upper bound in \eqref{Eq:Vio1}  below is  also obtained in \cite[Proposition 5.2]{CGLP20}. 
   
\begin{proposition} \label{prop3}
For an isotropic and mean square continuous  Gaussian random field $\{ Z(\x), \x \in \Md \}$ with mean 0 and 
covariance function (\ref{cov.mf1}), if (\ref{thm2.ineq2})  holds for a fixed $l_0 \in \N$ and positive constants 
$\gamma_1, \gamma_2$,  and $\nu \in (0, 2)$, then there are positive constants $\delta_0$, $K_1$ and $K_2$ 
such that
 \begin{equation}
 \label{Eq:Vio1}
 K_1 \rho (\x_1, \x_2)^\nu  \le \rE \left( Z(\x_1)- Z(\x_2) \right)^2  \le K_2 \rho (\x_1, \x_2)^\nu,
 \end{equation}
 holds   for all $\x_1, \x_2 \in  \Md$ with 
   $\rho (\x_1, \x_2)\le \delta_0$.
 \end{proposition}

 Proposition \ref{prop3} implies that many regularity and fractal properties of $\{ Z(\x), \x \in \Md \}$ are determined by 
 the index $\nu$. For example, the upper bound in (\ref{Eq:Vio1}) and  Kolmogorov's continuity theorem together  
 imply that $Z(\x)$ is H\"older continuous of any order $< \nu/2$ (see, e.g., \cite[Corollary 5.3]{CGLP20}). It can also 
 be  shown in a standard way that the Hausdorff dimension of the trajectory (the graph set) of $Z$,  Gr$Z(\Md) 
 = \{(\x, Z(\x)), \x \in \Md\} \subseteq \Md\times \R$, 
 is given by 
 \[
 \dim_{\rm H} \hbox{Gr} Z(\Md) = d + 1 -  \nu/2\, , \quad {\rm a.s.},
 \]
where $\dim_{\rm H} $ denotes Hausdorff dimension \cite{Fal04}. Because of 
these results, we call $\nu/2$ the fractal index of the random field  $\{ Z(\x), \x \in \Md \}$. When one uses such a 
random field as a statistics model to fit data sampled from values defined on $\Md$, it will be important to estimate 
the fractal index $\nu/2$. We refer to  \cite{MPbook11} for more information on statistical inference of random fields 
on the sphere and to \cite{Bhattacharya2012} for general nonparametric theory of statistics on manifolds.

{\it Remark:}
 If (\ref{thm2.ineq2}) holds for some $\nu > 2$, then it can be shown that the sample function of $\{ Z(\x), \x \in \Md \}$ is 
 continuously differentiable. We will not pursue this further in the present paper. 

 \begin{example}
For $\nu \in (0, 2]$, it is shown in Example 4 of \cite{LuMa2020} that
$$ C (\rho(\x_1,\x_2)) = 1-\left(\sin\frac{\rho(\x_1,\x_2)}{2} \right)^\nu, ~~~~ \x_1,\x_2\in\Md, $$
is the covariance function of an isotropic Gaussian random field on $\Md$. We will show that its spectral coefficients 
$b_n\sim n^{-\nu-1}$ on $\Md$. Denote the expansion of a bounded zonal function $f(x)$ on $\Md$ by
$$ f(x)=\sum_{n=0}^\infty b^{\alpha,\beta}_n(f)p^{\alpha,\beta}_n(\cos x),$$
where $p^{(\alpha,\beta)}_n(x)=P^{(\alpha,\beta)}_n(x)/P^{(\alpha,\beta)}_n(1)$.

Now we prove by induction that for $\nu>0$ and all integers $n \ge 0$,
\be\label{sinH}
b^{\alpha,\beta}_n(\sin^\nu\hf{x})=\frac{(2n+\alpha+\beta+1)\Gamma(\alpha+H+1)\Gamma(n+\alpha+\beta+1)\Gamma(n-H)}
{n!\Gamma(\alpha+1)\Gamma(\alpha+\beta+n+H+2)\Gamma(-H)}, 
\ee
where $H=\nu/2$. For $n=0$, $p^{\alpha,\beta}_0(x)=1$ implies
$$ b^{\alpha,\beta}_0(\sin^\nu\hf{x})=\frac{\Gamma(\alpha+H+1)\Gamma(\alpha+\beta+2)}{\Gamma(\alpha+1)\Gamma(\alpha+\beta+H+2)}.
$$
Thus \eq{sinH} holds for $n=0$. We assume that \eq{sinH} holds for some $n \ge 0$. Then,  by (18.9.6) in \cite{Olver2010},
$$ \frac{(n+1)(n+\beta+1)}{2n+\alpha+\beta+3}b^{\alpha,\beta}_{n+1}(f)=\frac{(n+\alpha+\beta+1)(n+\alpha+1)}
{2n+\alpha+\beta+1}b^{\alpha,\beta}_n(f)-(\alpha+1)b^{\alpha+1,\beta}_n(f), $$
and the induction assumption, we derive
\bes
b^{\alpha,\beta}_{n+1}(\sin^\nu\hf{x})&=&\frac{2n+\alpha+\beta+3}{(n+1)(n+\beta+1)}\frac{n+\alpha+\beta+1}
{2n+\alpha+\beta+1}\frac{(n+\beta+1)(n-H)}{\alpha+\beta+n+H+2}b^{\alpha,\beta}_n(\sin^\nu\hf{x})\non\\
&=&\frac{(2n+\alpha+\beta+3)\Gamma(\alpha+H+1)\Gamma(n+\alpha+\beta+2)\Gamma(n+1-H)}
{(n+1)!\Gamma(\alpha+1)\Gamma(\alpha+\beta+n+H+3)\Gamma(-H)}.\non
\ees
Therefore  \eq{sinH} holds for all integers $n\ge 0$. Consequently, %By \eq{sinH},
$$ b^{\alpha,\beta}_n(1-\sin^\nu\hf{x})=\delta_{n0}-\frac{(2n+\alpha+\beta+1)\Gamma(\alpha+H+1)
\Gamma(n+\alpha+\beta+1)\Gamma(n-H)}{n!\Gamma(\alpha+1)\Gamma(\alpha+\beta+n+H+2)\Gamma(-H)}, $$
which gives
$$ b^{\alpha,\beta}_n(1-\sin^\nu\hf{x})\sim \frac{-2\Gamma(\alpha+H+1)}{\Gamma(\alpha+1)\Gamma(-H)n^{\nu+1}}. $$
 \end{example}
 
 \begin{example}
Let
$$ p_0(x)=(\pi-x)^2,\ \ \ p_1(x)=2\pi^2(\pi-x)^2-(\pi-x)^4. $$
The functions $C(\xx_1,\xx_2)=p_0(\rho(\xx_1,\xx_2))$ and $C(\xx_1,\xx_2)=p_1(\rho(\xx_1,\xx_2))$ are 
covariance functions of isotropic Gaussian random fields on $\Md$, with spectral coefficients $b_n\sim n^{-2}$ 
and $b_n\sim n^{-4}$, respectively. For $\beta=-1/2$, integration by parts shows that
$$ b^{\alpha,-1/2}_n(p_0(x))=\frac{2\sqrt{\pi}\Gamma(n)}{n\Gamma(n+\half)}\frac{2n+\alpha+1}{n+\alpha+\half}
\frac{\Gamma(\alpha+\hf{3})\Gamma(n+\alpha+1)}{\Gamma(\alpha+1)\Gamma(n+\alpha+\hf{3})}, $$
and
$$ b_n^{\alpha,-\half}(p_1(x))=b_n^{\alpha,-\half}(p_0(x))12\left(\sum_{k=0}^\infty\inv{(n+k)^2}- 
\sum_{k=0}^\infty\inv{(n+k+\alpha+\hf{3})^2}\right). $$
Therefore,
$$ b^{\alpha,-1/2}_n(p_0(x))\sim \frac{4\sqrt{\pi}\Gamma(\alpha+\hf{3})}{\Gamma(\alpha+1)n^2},\ \ \ 
b^{\alpha,-1/2}_n(p_1(x))\sim \frac{48\sqrt{\pi}\Gamma(\alpha+\hf{5})}{\Gamma(\alpha+1)n^4}. $$
Following Example 1 of \cite{LuMa2020}, it can be shown that $p_0(x)$ and $p_1(x)$ are positive definite 
functions on all $\Md$, which implies that their asymptotic spectral coefficients are the same on all $\Md$ 
with the same dimension $d$.
 \end{example}

\section{Proofs}
 
In this section, we provide proofs for our main results, in the order of Propositions 1-3 and Theorems 1-2.

 \subsection{Proof of Proposition \ref{prop2} }
 
 Suppose that $C(\x_1, \x_2)$ is the covariance function of a Gaussian random field $\{ Z(\x),$ $ \x \in \D \}$. 
 For every $n \in \N$, any $\x_k \in \D$
 and any $a_k \in \R$ ($k = 1, \ldots, n$),  by applying the Cauchy-Schwarz inequality we obtain
 \[
   \begin{split} 
   &\left\{ \rE  \left( (Z (\x_0)- \rE Z(\x_0))  \sum_{k=1}^n a_k ( Z(\x_k) -\rE Z(\x_k)) \right) \right\}^2 \\
   &\le \var (Z(\x_0)) \var \left( \sum_{k=1}^n a_k ( Z(\x_k) -\rE Z(\x_k)) \right), 
         \end{split} 
\]         
or
$$  
\sum_{i=1}^n \sum_{j=1}^n a_i a_j   C( \x_i, \x_0) C (\x_j, \x_0 )  \le  C( \x_0, \x_0) \sum_{i=1}^n \sum_{j=1}^n a_i a_j C(\x_i, \x_j). 
$$
This implies that $ C(\x_0, \x_0) C(\x_1, \x_2) -C(\x_1, \x_0) C(\x_2, \x_0)$, $ \x_i, \x_2 \in \D$, is a positive definite 
function, and thus is a covariance function on $\D$.
 
 \subsection{Proof of Proposition  \ref{prop1} }
 
 To define a specific function on $[0, \pi]$ satisfying conditions in (i)-(iii),  let
\be
g_\eps(\vartheta)= 
  \left\{
   \begin{array}{ll}
    \phi\left(\frac{\sin\hf{\vartheta}}{\sin\hf{\eps}}\right),
     ~  &  ~  ~ \mbox{if} ~ \beta ~ \mbox{is an integer}, \\
      \cos \hf{\vartheta} \phi\left(\frac{\sin\hf{\vartheta}}{\sin\hf{\eps}}\right), ~ &
     ~ \mbox{if} ~ \beta- \frac{1}{2} ~ \mbox{is an integer},
\end{array}   \right.  
\ee
where 
   $$  \phi (x)  = \left\{ 
                  \begin{array}{ll}   
                   (1-x^2)_+^R\frac{P_{2K}(x)}{P_{2K}(0)},       ~ &  ~ \mbox{if} ~ d ~ \mbox{is odd}, \\  
                   (1-x^2)_+^R P_K(1-2x^2),      ~ &  ~ \mbox{if} ~ d ~ \mbox{is even},
                    \end{array}   \right. $$
in which $P_K (x) $ is the Legendre polynomial of degree $K$, with integers $R$ and $K$ such that
$$ R\ge r+\alpha-\half,\ \ \ K=n_0+R+\lceil\alpha\rceil+\lceil\beta\rceil. $$
It's clear that $g_\eps (\vartheta)$ is continuous on $[0, \pi]$, $g_\eps(0)=1$, and $g_\eps(\pi)=0$. 

For an odd $d$,  $\beta+\frac{1}{2}$ is an integer,  as is seen from Table 1. Making the transform 
$y=\sin\hf{\vartheta}/\sin\hf{\eps}$ 
we obtain
 \begin{eqnarray*}
 &  & \int_0^\pi g_\eps(\vartheta)P^{ (\alpha, \beta)}_n(\cos\theta)\sin^{2\alpha+1}\hf{\theta}
 \cos^{2\beta+1}\hf{\vartheta} d \vartheta \\
 &  =  &  \int_0^\eps    \cos \hf{\vartheta}   
           \left( 1- \frac{\sin^2 \hf{\vartheta}}{\sin^2\hf{\eps}} \right)^R
             \frac{P_{2K} \left( \frac{\sin \hf{\vartheta}}{\sin \hf{\eps}}   \right)}{P_{2K} (0)}
               P^{ (\alpha, \beta)}_n(\cos\theta)\sin^{2\alpha+1}\hf{\theta}\cos^{2\beta+1}\hf{\theta} d \vartheta \\
  & = &  \frac{ 2 \sin^{2 \alpha+2} \hf{\eps} }{P_{2K} (0)} 
              \int_0^1  (1-y^2)^R  y^{ 2 \alpha+1} \left( 1- y^2 \sin^2 \hf{\eps} \right)^{ \beta+\frac{1}{2}  } P_{2 K}(y)  
              P^{ (\alpha, \beta)}_n \left(    1- 2y^2 \sin^2 \hf{\eps}   \right)   dy \\
 & = & \int_{-1}^1 P_{2K}(y) h_1(y^2)dy, 
 \end{eqnarray*}
 where  
   $$ h_1 (y^2) = \frac{ \sin^{2 \alpha+2} \hf{\eps} }{P_{2K} (0)} 
                (1-y^2)^R y^{ 2 \alpha+1} \left( 1- y^2 \sin^2 \hf{\eps} \right)^{ \beta+\frac{1}{2}  } 
              P^{ (\alpha, \beta)}_n \left(    1- 2y^2 \sin^2 \hf{\eps}   \right) $$
              is  a polynomial of $y$ of degree $2n+2R+ 2 \alpha +2 \beta+1$.
  For $0 \le n  < n_0$,  we have 
  $$2 K =  2 n_0+2R+ 2 \alpha+2 \beta+1 > 2n+2R+2 \alpha +2 \beta+ 1, $$ 
  so that $\int_{-1}^1 P_{2K}(y) h_1 (y^2)dy =0$ by the orthogonality of the Legendre polynomial. 
  It implies that 
  $b^{(\alpha,\beta)}_n(g_\eps)=0$ for $0 \le n< n_0$.           
              
For an even $d$ and $\beta \neq -\frac{1}{2}$,  $\alpha$ and $\beta$ are integers, as is seen from 
Table \ref{table1}.  By making the transform  $y=\sin\hf{\vartheta}/\sin\hf{\eps}$ and followed by a change 
of variable $w= 1-2 y^2$,  we obtain
 \begin{eqnarray*}
 &  & \int_0^\pi g_\eps(\vartheta)P^{ (\alpha, \beta)}_n(\cos\theta)\sin^{2\alpha+1}\hf{\theta}\cos^{2\beta+1}
 \hf{\vartheta} d \vartheta \\
 &  =  &  \int_0^\eps      
           \left( 1- \frac{\sin^2 \hf{\vartheta}}{\sin^2\hf{\eps}} \right)^R
             P_{K} \left( 1- \frac{ 2 \sin^2 \hf{\vartheta}}{\sin^2 \hf{\eps}}   \right)
               P^{ (\alpha, \beta)}_n(\cos\theta)\sin^{2\alpha+1}\hf{\theta}\cos^{2\beta+1}\hf{\theta} d \vartheta \\
  & = &  2  \sin^{2 \alpha+2} \hf{\eps}
              \int_0^1  (1-y^2)^R y^{ 2 \alpha+1} \left( 1- y^2 \sin^2 \hf{\eps} \right)^{ \beta  } P_{ K} \left( 1-2 y^2 \right)  
              P^{ (\alpha, \beta)}_n \left(    1- 2y^2 \sin^2 \hf{\eps}   \right)   dy \\
 & = &   \frac{1}{2}  \sin^{2 \alpha+2} \hf{\eps}
              \int_{-1}^1  \left( \frac{1+w}{2} \right)^R  \left(\frac{1-w}{2} \right)^{  \alpha} \left( 1-  \frac{1-w}{2} 
              \sin^2 \hf{\eps} \right)^{ \beta  }
              P_K(w) P^{ (\alpha, \beta)}_n \left(    1- (1-w) \sin^2 \hf{\eps}   \right)   d w \\
  & = &    \int_{-1}^1 P_K(w) h_2 (w) dw,          
 \end{eqnarray*}
where $h_2 (w) $ is the polynomial of degree of $n+R +\alpha+\beta$ defined by
\[
\begin{split}
h_2 (w) &= \frac{1}{2}  \sin^{2 \alpha+2} \hf{\eps}
\left( \frac{1+w}{2} \right)^R  \left(\frac{1-w}{2} \right)^{  \alpha} \left( 1-  \frac{1-w}{2} \sin^2 \hf{\eps} \right)^{ \beta} \\
&\qquad \qquad \qquad \times P^{ (\alpha, \beta)}_n \left(    1- (1-w) \sin^2 \hf{\eps}   \right). 
\end{split}
\]
         
For $0 \le n < n_0$,\, it follows from  the orthogonality of the Legendre polynomial and $K=n_0+R+\alpha+\beta > 
n+R +\alpha+\beta$ that $ \int_{-1}^1 P_K(w) h_2 (w) dw =0$, and  $b^{ (\alpha,\beta)}_n(g_\eps)=0$. 

Now we prove $|b^{ (\alpha,\beta)}_n(f_\eps)|\le M_r\eps(1+n\eps)^{-r}$. By  Theorem 8.1.1 of  \cite{Szego1975},  
$$ \lim_{n\to\infty}\left(\frac{x}{2n}\right)^\alpha P^{(\alpha,\beta)}_n \left( \cos \frac{x}{n} \right)=J_\alpha(x), $$
where $J_\alpha(x)$ is the Bessel function of the first kind, we have
\be\label{clim}
\lim_{\eps\to 0,\ n\eps=k}\frac{b^{(\alpha,\beta)}_n(g_\eps)}{\eps}
=\frac{2k^{2\alpha+1}}{\Gamma(\alpha+1)}I^\alpha_k(\phi),
\ee
where
$$ I^\alpha_k(\phi)=\int_0^1\phi(y)\frac{J_\alpha(ky)}{k^\alpha}\left(\hf{y}\right)^{\alpha+1}dy.$$
In particular, for $\alpha=- \frac{1}{2}$,
$$ I^{-\half}_k(\phi)=\inv{\sqrt{\pi}}\int_0^1\phi(y)\cos(ky)dy. $$
Since $\phi (y) $ is an even polynomial with derivatives up to order $R-1$ vanishing at 1, integration by parts gives
$$ I^{-\half}_k(\phi)=\inv{\sqrt{\pi}k^{R+1}}\left(\phi^{(R)}(1)g_R(k)-\int_0^1\phi^{(R+1)}(y)g_R(ky)dy\right), $$
where
   $$ g_R (y) = \left\{
                        \begin{array}{ll}
                         (-1)^{\lfloor R/2\rfloor}\sin(y), ~  &  ~ \mbox{if} ~ R ~ \mbox{is even}, \\
                         (-1)^{\lfloor R/2\rfloor}\cos(y), ~  &  ~ \mbox{if} ~ R ~ \mbox{is odd}.
                         \end{array}    \right. $$
     For $\alpha=0$, integrating by parts and using the derivative formulas,
$$  (J_1(x)x)'=J_0(x)x,\ \ \ J_0(x)'=-J_1(x),  $$
we get
  \begin{eqnarray*}
      I^0_k(\phi) & = & \half\int_0^1\phi(y)J_0(ky)ydy  \\
           &  =  &  \inv{2k^{R+1}}\left(\phi_R(1)h_R(k)-\int_0^1\phi_{R+1}(y)h_R(ky)ydy\right),
  \end{eqnarray*}
where 
   $$ h_R (y) = \left\{
                        \begin{array}{ll}
                         (-1)^{\lfloor R/2\rfloor} J_1(y), ~  &  ~ \mbox{if} ~ R ~ \mbox{is even}, \\
                         (-1)^{\lfloor R/2\rfloor} J_0(y), ~  &  ~ \mbox{if} ~ R ~ \mbox{is odd},
                         \end{array}    \right. $$
                          and $\phi_m(y)$ are defined recursively by $\phi_0(y)=\phi(y)$, and for $m\ge0$,
  $$ \phi_{2m+1}(y)=\phi_{2m}'(y),\ \ \ \phi_{2m+2}(y)=(y\phi_{2m+1}(y))'/y,  $$
which are well-defined since $\phi (y) $ is an even polynomial. For higher dimensions, using the 
derivative formula $(J_\alpha(x)/x^\alpha)'=-J_{\alpha+1}(x)/x^\alpha$, we get
      $$     I^{\alpha+1}_k(\phi)=-\frac{d}{2kdk}I^\alpha_k(\phi), $$
on odd-dimensional spaces, for large $k$,
     $$  I^\alpha_k(\phi)\sim\frac{\phi^{(R)}(1)g_{R+\lceil\alpha\rceil}(k)}{\sqrt{\pi}k^{R+1}(-2k)^{\lceil\alpha\rceil}}, $$
and 
 on even-dimensional spaces, for large $k$,
     $$   I^\alpha_k(\phi)\sim\frac{\phi_R(1)h_{R+\lceil\alpha\rceil}(k)}{2k^{R+1}(-2k)^{\lceil\alpha\rceil}}.  $$
Since the derivatives of $\phi$ up to order $R-1$ vanishes at 1, $\phi_R(1)=\phi^{(R)}(1)$. Using the asymptotic 
form of Bessel functions,
   $$  J_\alpha(z)=\sqrt{\frac{2}{\pi z}}\left(\cos(z-\hf{\alpha}\pi-\frac{\pi}{4})+O \left( \inv{z} \right) \right), $$
we get the unified asymptotic form of $I^\alpha_k(\phi)$,
   $$  I^\alpha_k(\phi)=\frac{\phi^{(R)}(1)}{\sqrt{2\pi}2^\alpha k^{R+\alpha+\hf{3}}}
       \left(\sin(k+\hf{\pi}(R-\alpha-\half))+O \left( \inv{k} \right) \right), $$
and substituting it into \eq{clim} we get
    $$  b^{ (\alpha,\beta)}_n(g_\eps)=\eps\cdot O((n\eps)^{\alpha-R-\half})=\eps\cdot O((n\eps)^{-r}). $$ 
On the other hand, by \eq{cab},
  $$  b^{(\alpha,\beta)}_n(g_\eps)=O(n^{2\alpha+1}\eps^{2\alpha+2})=\eps\cdot O((n\eps)^{2\alpha+1}).  $$
Therefore $| b^{(\alpha,\beta)}_n(g_\eps)| \le M_r\eps(1+n\eps)^{-r}$.

\subsection{Proof of Proposition  \ref{prop3} }
  
To show (\ref{Eq:Vio1}), we  apply the following recurrence relation of the Jacobi polynomials \cite{Szego1975},
$$  p^{(\alpha,\beta)}_l(\cos\theta)-p^{(\alpha,\beta)}_{l+1}(\cos\theta)
=\frac{2l+\alpha+\beta+2}{\alpha+1}p^{(\alpha+1,\beta)}_l(\cos\theta)\sin^2\frac{\theta}{2}, $$
where $p^{(\alpha,\beta)}_l(x)=P^{(\alpha,\beta)}_l(x)/P^{(\alpha,\beta)}_l(1)$.
It   is known that 
$|P^{(\alpha,\beta)}_l(x)| \le P^{(\alpha,\beta)}_l(1)$, and 
$ 0 \le 1-p^{(\alpha,\beta)}_l(x) \le 2$ for $x \in [-1, 1]$. Also, $0 \le \frac{ \sin \frac{\theta}{2}}{\theta} \le \frac{\pi}{2}, \theta \in  (0, \pi]$. 
 It gives that
\begin{eqnarray*}
1-p^{(\alpha,\beta)}_l(\cos\theta)
 & = &   \sum_{j=0}^{l-1} (  p^{(\alpha,\beta)}_j(\cos\theta)-p^{(\alpha,\beta)}_{j+1}(\cos\theta) ) \\
 &   =  & \sum_{j=0}^{l-1}\frac{2m+\alpha+\beta+2}{\alpha+1}
p^{(\alpha+1,\beta)}_j(\cos\theta)\sin^2\frac{\theta}{2}\\
 & = & \frac{l(l-1)+\alpha+\beta+2}{\alpha+1} \sin^2\frac{\theta}{2}\\ 
& \le & \sum_{m=0}^{l-1}\frac{2j+\alpha+\beta+2}{\alpha+1}\sin^2\frac{\theta}{2} \le K\,l^2\theta^2, ~~~ 0 \le \theta \le \pi,
\end{eqnarray*}
for some constant $K>0$ that depends only on $\alpha$ and $\beta$. 
Hence, for $\theta=\rho(\x_1, \x_2)$,  by (\ref{Eq:viog1}) and the condition that $b_ll^{1+\nu}\le \gamma_2$, 
we obtain 
\be \label{Eq:Vario1}
\begin{split}
\rE \big[(Z(\x_1)- Z(\x_2))^2\big] &= 2\sum_{l=0}^\infty b_l \Big(1-p^{(\alpha,\beta)}_l(\cos\theta)\Big)\\
&\le 2\sum_{0\le l<1/\theta} b_lK_1l^2\theta^2+4\sum_{l\ge 1/\theta} b_l\\
&\le K_2\theta^\nu,
\end{split}
\ee
for some finite constant $K_2>0$. This proves the upper bound in (\ref{Eq:Vio1}). 

On the other hand,  the lower bound in  (\ref{Eq:Vio1}) follows from Theorem 1 because 
$$\rE \big[(Z(\x_1)- Z(\x_2))^2\big] \ge  \var ( Z(\x_1) |  Z(\x_2)) \ge \gamma \rho(\x_1, \x_2)^\nu.$$ 
Alternatively, it can be proved in the following elementary way. For $x\in[0,\pi]$,
$$ p^{\alpha,\beta}_n(\cos x)=\sum_{m=0}^n (-1)^m{n\choose m}\frac{\Gamma(n+\alpha+\beta+m+1)\Gamma(\alpha+1)}
{\Gamma(n+\alpha+\beta+1)\Gamma(\alpha+m+1)}\Big(\sin\hf{x}\Big)^{2m}. $$
Let
\bes
&\ &g(\cos x)=\frac{1-p^{\alpha,\beta}_n(\cos x)}{(\sin\hf{x})^2}\non\\
& \ & =\frac{n(n+\alpha+\beta+1)}{\alpha+1}-\sum_{m=2}^n (-1)^m{n\choose m}\frac{\Gamma(n+\alpha+\beta+m+1)
\Gamma(\alpha+1)}{\Gamma(n+\alpha+\beta+1)\Gamma(\alpha+m+1)}\Big(\sin\hf{x}\Big)^{2m-2}.\non
\ees
We have
\[
\begin{split}
& {n\choose m}\frac{\Gamma(n+\alpha+\beta+m+1)\Gamma(\alpha+1)}
{\Gamma(n+\alpha+\beta+1)\Gamma(\alpha+m+1)}\Big(\sin\hf{x}\Big)^{2m-2}\\
&\le \frac{(n+\alpha+\beta+1)^{2m}\Gamma(\alpha+1)}{m!\Gamma(\alpha+m+1)}\Big(\hf{x}\Big)^{2m-2}. 
\end{split}
\]
For any $\eps>0$, there exists $\delta>0$ such that for $n\ge\alpha+\beta+1$ and $0\le nx\le\delta$,
$$ 
\sum_{m=2}^n \frac{(n+\alpha+\beta+1)^{2m}\Gamma(\alpha+1)}{m!\Gamma(\alpha+m+1)}\Big(\hf{x}\Big)^{2m-2}
\le\sum_{m=2}^n \frac{4n^2(nx)^{2m-2}\Gamma(\alpha+1)}{m!\Gamma(\alpha+m+1)}\le n^2\eps. $$
Therefore, there exists $K_1>0$ for which
$$ 1-p^{\alpha,\beta}_n(\cos x)=g(\cos x) \Big(\sin\hf{x}\Big)^2\ge K_1(nx)^2. $$
Given that $b_n\ge n^{-\nu-1}$ for large $n$, there exists $K_2>0$ such that
$$ \sum_{n=0}^\infty b_n \left(1-p^{\alpha,\beta}_n(\cos x) \right)\ge \sum_{n=\alpha+\beta+1}^{\delta/x}\frac{K_1(nx)^2}{n^{\nu+1}}\ge K_2 x^\nu. $$

 \subsection{Proof of Theorem 1}
 
Write $\varepsilon = \min\limits_{1 \le  k \le n} \rho (\x, \x_k).$ 
For the  Gaussian random field $\{ Z(\x), \x \in  \Md \}$, in order to verify inequality \eqref{thm1.ineq2}, it suffices 
to show that there is a positive constant $\gamma$ such that
\begin{equation}
\label{inf.ineq}
 \rE \bigg( Z(\x) - \sum_{k=1}^n a_k Z(\x_k) \bigg)^2  \ge \gamma \varepsilon^\nu,
\end{equation}
 holds for all  $n \in \N$,  $ \x, \,\x_k \in \Md$  ($k=1, \ldots, n$) with $\min\limits_{1 \le  k \le n} \rho (\x, \x_k) > 0$ 
 and all $a_k \in \R$ ($k=1, \ldots, n$).
    We have
        \begin{eqnarray*}
        &  &    \rE \bigg( Z(\x) - \sum_{k=1}^n a_k Z(\x_k) \bigg)^2  \\
        & = &  C(\x, \x) - 2\sum_{i=1}^n a_i C(\x, \x_i)+\sum_{i=1}^n \sum_{j=1}^n a_i a_j  C(\x_i, \x_j)  \\
           & = & \sum_{l=0}^\infty b_l \left[  \left( 1-  \sum_{k=1}^n a_k
           \frac{P_l^{ \left( \alpha, \beta \right)} (\cos \rho  (\x, \x_k))}{  P_l^{ \left( \alpha, \beta  \right)} (1)}  \right)^2\right.\\
          & & + \left. \sum_{i=1}^n \sum_{j=1}^n a_i a_j
            \left( \frac{P_l^{ \left( \alpha, \beta  \right)} (\cos \rho  (\x_i, \x_j))}{  P_l^{ \left( \alpha, \beta \right)} (1)} -
             \frac{P_l^{ \left( \alpha, \beta \right)} (\cos \rho (\x, \x_i))}{  P_l^{ \left( \alpha, \beta \right)} (1)} 
              \frac{P_l^{ \left( \alpha, \beta \right)} (\cos \rho (\x, \x_j))}{  P_l^{ \left( \alpha, \beta \right)} (1)}  \right) \right]\\
          & \ge &   \sum_{l=0}^\infty b_l   \left( 
                       1-  \sum_{k=1}^n a_k
                          \frac{P_l^{ \left( \alpha, \beta \right)} (\cos \rho (\x, \x_k))}{  P_l^{ \left( \alpha, \beta \right)} (1)}  \right)^2, 
      \end{eqnarray*}
      where the last inequality holds since, for every $l \in \N$,
 \begin{eqnarray*}     
& & \sum_{i=1}^n \sum_{j=1}^n a_i a_j \left( 
  \frac{P_l^{ \left( \alpha, \beta  \right)} (\cos \rho (\x_i, \x_j))}{  P_l^{ \left( \alpha, \beta \right)} (1)}  -
    \frac{P_l^{ \left( \alpha, \beta \right)} (\cos \rho (\x, \x_i))}{  P_l^{ \left( \alpha, \beta  \right)} (1)}  
     \frac{P_l^{ \left( \alpha, \beta \right)} (\cos \rho (\x, \x_j))}{  P_l^{ \left( \alpha, \beta \right)} (1)} \right) \ge 0,
\end{eqnarray*} 
which is due to Proposition \ref{prop2}, while $\frac{P_l^{ \left( \alpha, \beta \right)} (\cos \rho (\x_1, \x_2))}
{  P_l^{ \left( \alpha, \beta \right)} (1)} $
is known to be a covariance function on $\Md$ by Lemma 3 of \cite{MaMalyarenko2018}.
     
 For a continuous function $g_\eps (\vartheta)$ on $[0, \pi]$ satisfying conditions (i)-(iii) in Proposition \ref{prop1}  
 with $r >1  + \nu/2$, we consider
 $$ I =    \sum_{l=0}^\infty b_l^{(\alpha, \beta)}  (g_\eps)  \left(1
         - \sum_{k=1}^n  a_k   \frac{P_l^{ \left( \alpha, \beta \right)} (\cos \rho  (\x, \x_k))}{  P_l^{ \left( \alpha, \beta \right)} (1)} \right). 
 $$
         On one hand,  it follows from $\rho  (\x, \x_k) \ge \eps$   and Proposition \ref{prop1} (i) that
         $g_\eps (\rho (\x, \x_k))  =0$ ($k = 1, \ldots, n$), so that 
  \begin{equation*}
  \begin{split}
    I  & =    \sum_{l=0}^\infty    b_l^{(\alpha, \beta)}  (g_\eps)
                  - \sum_{k=1}^n  a_k  \sum_{l=0}^\infty  b_l^{(\alpha, \beta)}  (g_\eps)
                   \frac{P_l^{ \left( \alpha, \beta \right)} (\cos \rho  (\x, \x_k))}{  P_l^{ \left( \alpha, \beta \right)} (1)} \\
       & =    g_\eps (0) -  \sum_{k=1}^n  a_k  g_\eps (\rho (\x, \x_k))  \\
       & =     1.
       \end{split}
   \end{equation*}   
On the other hand, an application of the Cauchy-Schwarz inequality yields that
      \begin{eqnarray*}
      I^2  & = &  \left\{  \sum_{l=l_0}^\infty  \frac{b_l^{(\alpha, \beta)}  (g_\eps)}{ \sqrt{b_l} }   \sqrt{b_l} \left( 1
         - \sum_{k=1}^n  a_k   \frac{P_l^{ \left( \frac{d-1}{2} \right)} (\cos \rho (\x, \x_k))}{  P_l^{ \left( \alpha, \beta \right)} (1)}  \right)    \right\}^2  \\
         & \le &   \sum_{l=l_0}^\infty   \frac{ \left( b_l^{(\alpha, \beta)}  (g_\eps)  \right)^2}{ b_l}
                     \sum_{l=l_0}^\infty  b_l \left( 1
         - \sum_{k=1}^n  a_k   \frac{P_l^{ \left( \alpha, \beta \right)} (\cos \rho (\x, \x_k))}{  P_l^{ \left( \alpha, \beta \right)} (1)}  \right)^2 \\
          & \le & \frac{1}{M_r \eps^\nu} \sum_{l=l_0}^\infty  b_l \left( 1
         - \sum_{k=1}^n  a_k   \frac{P_l^{ \left( \frac{d-1}{2} \right)} (\cos \vartheta (\x, \x_k))}{  P_l^{ \left( \frac{d-1}{2} \right)} (1)}  \right)^2,
      \end{eqnarray*}   
 where the last inequality follows from $|b^{(\alpha,\beta)}_l (g_\eps)|\le M_r \eps (1+ l \eps)^{-r}$ by Proposition \ref{prop1} (iii),  
 $b_l  (1+l)^{1+\nu}   \ge K$ for $l \ge l_0$ by   (\ref{thm1.ineq1}), and
 \[
 \begin{split}
   \sum_{l=l_0}^\infty  \frac{ \left( b_l^{(\alpha, \beta)}  (g_\eps)  \right)^2}{ b_l} 
       &\le  \sum_{l=0}^\infty     \frac{ M_r^2  \eps^2  (1+ l \eps)^{-2 r} }{ K (1+l)^{ -1-\nu}  } \\ 
       &\le\frac{M_r^2 }{K} \eps^2  \int_0^\infty \frac{(1+x)^{1+\nu}}{(1+\eps x)^{2r}}\, dx
       = \frac{M_r^2}{K}\frac{\eps^{-\nu}}{2r-\nu-2}.
       \end{split}
       \]       
Consequently, inequality (\ref{thm1.ineq2}) is obtained.

 \subsection{Proof of Theorem 2}
 
We start with the following zero-one law, which is proved by applying the Karhunen-Lo\`eve expansion
for $Z$ (cf. \cite[Chapter 2]{Malyarenko2013}) and the Kolmogorov's zero-one law. 

\begin{lemma}\label{Lem:01law}
Let $Z= \{Z(\x), \, \x\in \Md\}$ be a centered isotropic Gaussian random field on $\Md$ that satisfies the 
conditions of Theorem \ref{thm2}.
Then there is a constant $K \in [0, \infty]$ such that 
\begin{equation}\label{Eq:UM0}
\lim_{\varepsilon \to 0} \sup_{^{\;\; \x_1,\x_2\in {\Md}}_{\rho(\x_1, \x_2)\leq\varepsilon} } 
\frac{|Z(\x_1)-Z(\x_2)|} {  \rho (\x_1, \x_2)^{\nu/2} \sqrt{ \ln  \rho (\x_1, \x_2)} } =K,\;\; \hbox{ a.s.}  
\end{equation}
\end{lemma}
\begin{proof}

Recall from Malyarenko 
\cite[Chapter 2]{Malyarenko2013} that $ \{ Z(\x), \x \in \Md \}$ has the following Karhunen-Lo\`eve expansion
\begin{equation}\label{Eq:KL}
Z(\x) = C_d \sum_{l = 0}^\infty \sum_{m=1}^{h(\Md,l)} \sqrt{\frac{b_l}{h(\Md, l)}} \ X_{l, m} Y_{l, m}(\x),
\end{equation}
with convergence in $L^2(\Omega,  L^2(\Md))$. In the above, $C_d = \sqrt{\frac{2 \pi^{(d+1)/2}} {\Gamma((d+1)/2)} }$ 
if $\Md = \mathbb{S}^d$ and $C_d = 1$ in all other cases in Figure 1; $h(\Md,l)$ is given by 
\[
h(\Md,l) = \frac{(2 l + \alpha + \beta + 1)\Gamma(\beta+1)\Gamma(l + \alpha + \beta + 1) \Gamma(l + \alpha  + 1)}
{\Gamma(\alpha+1)\Gamma(\alpha +\beta+2)l! \Gamma(l+ \beta+1)},
\]
where $\alpha$ and $\beta$ are the parameters given by Figure 1,  $\{b_l, l \ge 0\}$ is the  angular power spectrum 
of $Z$, and the sequences $\{X_{l, m}\}$ , and $\{Y_{l m}\}$ are specified as follows:
\begin{itemize}
 \item $\{X_{l, m}\}$ is a sequence of i.i.d. standard normal random variables.
 \item  $\{Y_{l ,m}\}$ are the eigenfunctions of the Laplace-Beltrami operator $\Delta_{\Md}$ on $\Md$, i.e.,
 \[
 - \Delta_{\Md}  Y_{l ,m} = \lambda_l Y_{l ,m},
 \]
where the eigenvalues $\lambda_l = \lambda_l (\alpha, \beta) = l(l + \alpha + \beta + 1), $ for all $l \in \mathbb{N}_0$.  
\end{itemize}

For every $l \ge 0$, the eigenfunctions $\{Y_{l ,m}, 1 \le m \le h(\Md,l)\}$ corresponding to the same
eigenvalue $\lambda_l$  form a finite dimensional vector space of dimension $h(\Md,l)$.
It is known that, for every $(l, m)$, the eigenfunction $Y_{l, m}(\x)$ is continuously differentiable and $\Md$ is compact. Recall that $\nu \in (0, 2)$. Hence for every  integer $L \ge 0$, 

\begin{equation*} 
\lim_{\varepsilon \to 0} \sup_{^{\;\; \x_1,\x_2\in {\Md}}_{\rho(\x_1, \x_2)\leq\varepsilon} } 
\frac{|Z_L(\x_1)-Z_L(\x_2)|} {  \rho (\x_1, \x_2)^{\nu/2} \sqrt{ \ln  \rho (\x_1, \x_2)} } =0,\;\; \hbox{ a.s.}, 
\end{equation*}
where
 \[
  Z_L(\x) =   C_d \sum_{l = 0}^L \sum_{m=1}^{h(\Md,l)} \sqrt{\frac{b_l}{h(\Md, l)}} \ X_{l, m} Y_{l, m}(\x).
  \]
Hence for every constant $\kappa_1\ge 0$, the event 
\[
E_{\kappa_1} = \Bigg\{\lim_{\varepsilon \to 0} \sup_{^{\;\; \x_1,\x_2\in {\Md}}_{\rho(\x_1, \x_2)
\leq\varepsilon} } 
\frac{|Z(\x_1)-Z(\x_2)|} {  \rho (\x_1, \x_2)^{\nu/2} \sqrt{ \ln  \rho (\x_1, \x_2)} } \le \kappa_1\Bigg\}
\]
is a tail event with respect to $\{X_{l, m}\}$. By Kolmogorov's zero-one law, we have $\P(E_{\kappa_1})$ $ = 0$ or 1. 
This implies \eqref{Eq:UM0} with $K = \sup\{\kappa_1\ge 0: \P(E_{\kappa_1}) = 0\}$.
\end{proof}
 
Now we prove Theorem 2. 

{\it Proof of Theorem 2}\ Because of the zero-one law in Lemma \ref{Lem:01law}, it is sufficient to prove the 
existence of positive and finite constants 
 $K_{5}$  and $K_{6}$  such that 
 \begin{equation}
 \label{Eq:UM1}
  \lim_{\varepsilon \to 0} \sup\limits_{ \substack{ \x_1, \x_2 \in \Md  \\  \rho (\x_1, \x_2) \le \varepsilon } }
    \frac{ | Z(\x_1) -Z(\x_2) |}{  \rho (\x_1, \x_2)^{\nu/2} \sqrt{ \ln  \rho (\x_1, \x_2)} } \le K_5\ \ \ \hbox{a.s.}  
 \end{equation}
 and
 \begin{equation}
 \label{Eq:UM2}
  \lim_{\varepsilon \to 0} \sup\limits_{ \substack{ \x_1, \x_2 \in \Md  \\ \rho (\x_1, \x_2) \le \varepsilon } }
   \frac{ | Z(\x_1) -Z(\x_2) |}{  \rho (\x_1, \x_2)^{\nu/2} \sqrt{ \ln  \rho (\x_1, \x_2)} } \ge K_6 \ \  \ \hbox{a.s.} 
 \end{equation}

The proof of  (\ref{Eq:UM1}) is quite standard. By Proposition \ref{prop3}, the canonical metric $d_Z$ defined by 
$$ d_Z(\x_1, \x_2) =  \sqrt{\rE \big[(Z(\x_1)- Z(\x_2))^2\big] }$$
satisfies $ d_Z (\x_1, \x_2) \le K\,\rho(x_1, \x_2)^{\nu/2}$ for all $\x_1, \x_2 \in \Md$. This implies that for every $\varepsilon \in (0, \pi)$, we have 
\begin{equation}\label{Eq:ent}
N(\Md, d_Z, \varepsilon) \le K \varepsilon^{-\frac{2 d} \nu},
\end{equation}
where $N(\Md, d_Z, \varepsilon) $ denotes the minimum number of $d_Z$-balls of radius $\varepsilon$ that are needed 
to cover $\Md$. Hence,  \eqref{Eq:UM1} follows from \eqref{Eq:ent} and Theorem 1.3.5 in \cite{RFG}.
 
For any $n\geq n_0$,  we choose a sequence of $2^{n}$ points $\{x_{n,i},1\leq
i \leq 2^{n}\} \subseteq \Md$  that are (approximately) equally separated in the
following sense: For every $2\leq k\leq 2^{n}$, we have 
\begin{equation} \label{Eq: xpoints}
\begin{split}
&\min_{1\leq i\leq k-1} \rho (\x_{n,k},\, \x_{n,i}) =
\rho(\x_{n,k}, \x_{n,k-1}) \\
 &\ \ K' 2^{-n} \le \rho(\x_{n,k}, \x_{n,k-1}) \le K 2^{-n}.
 \end{split} 
\end{equation}

With the choice of $\{\x_{n,i},1\leq i\leq 2^{n}\}$, we now prove  (\ref{Eq:UM2}) in a way 
that is similar to the proof in \cite{Lan2018}.
Notice that 
\begin{equation} \label{Eq:Umod-lower2}
\begin{split}
\lim_{\varepsilon \rightarrow 0} \sup_{\substack{ \x,\y\in \Md,  \\ 
\rho(\x,\y)\leq \varepsilon }}\frac{|Z(\x)-Z(\y)|}{\rho(\x,\y)^{\nu/2}\sqrt{|\ln \rho (\x,\y)|}}
 \geq \underset{n\rightarrow \infty }{\, \lim \inf }\max_{2\leq k\leq 2^{n}}
\frac{|Z(\x_{n,k})-Z(\x_{n,k-1})|}{2^{-n\nu/2}\sqrt{n}}.
\end{split}
\end{equation}
It is sufficient to prove that, almost surely, the last limit in (\ref{Eq:Umod-lower2}) is bounded 
below by a positive constant. This is done by applying the property of strong local nondeterminism 
in Theorem 1 and a standard Borel-Cantelli argument.
  
 Let $\eta >0$ be a constant whose value will be chosen later. We consider
the events 
\begin{equation*}
A_{m}=\bigg\{\max_{2\leq k\leq m}\big |Z(\x_{n,k})-Z(\x_{n,k-1})\big |\leq
\eta 2^{-n\nu/2}\sqrt{n}\bigg\}
\end{equation*}
for $m=2,3,\ldots ,2^{n}$. By conditioning on $A_{2^{n}-1}$ first, we can
write 
\begin{equation} \label{Eq:inter-1}
\begin{split}
{\mathbb{P}}\big(A_{2^{n}}\big)& ={\mathbb{P}}\big(A_{2^{n}-1}\big) \\
& \quad \times {\mathbb{P}}\bigg\{\big |Z(\x_{n,2^{n}})-Z(\x_{n,2^{n}-1})\big |
\leq \eta 2^{-n\nu/2}\sqrt{n}\, \big |A_{2^{n}-1}\bigg\}.
\end{split}
\end{equation}
 
Recall that, given the random variables in $A_{2^{n}-1}$, the conditional
distribution of the Gaussian random variable $
Z(\x_{n,2^{n}})-Z(\x_{n,2^{n}-1}) $ is still Gaussian, with the corresponding
conditional mean and variance as its mean and variance. By Theorem 1 and
\eqref {Eq: xpoints}, the aforementioned conditional variance satisfies 
\begin{equation*}
\mathrm{Var}\big(Z(\x_{n,2^{n}})-Z(\x_{n,2^{n}-1})\big|A_{2^{n}-1}\big)\geq
\gamma_3\,2^{-n \nu}.
\end{equation*}%
This and Anderson's inequality (see \cite{A55}) imply 
\begin{equation} \label{Eq:inter-2}
\begin{split}
& {\mathbb{P}}\bigg\{\big |Z(\x_{n,2^{n}})-Z(\x_{n,2^{n}-1})\big |\leq \eta
2^{-n\nu/2}\sqrt{n}\big |\,A_{2^{n}-1}\bigg\} \\
& \leq {\mathbb{P}}\bigg\{N(0,1)\leq \frac{\eta}{\sqrt{\gamma_3}}\, \sqrt{n}\bigg\} \\
& \leq 1-\frac{\sqrt{\gamma_3}}{\eta \sqrt{n}}\exp \bigg(-\frac{\eta ^{2}n}{2\gamma_3}\bigg) \\
& \leq \exp \bigg(-\frac{\sqrt{\gamma_3}}{\eta \sqrt{n}}\exp \Big(-\frac{\eta ^{2}n}{%
2 \gamma_3}\Big)\bigg).
\end{split}
\end{equation}%
In deriving the last two inequalities, we have applied Mill's ratio and the
elementary inequality $1-x\leq e^{-x}$ for $x>0$. Iterating this procedure
in (\ref{Eq:inter-1}) and \eqref{Eq:inter-2} for $2^{n}-1$ more times, we
obtain 
\begin{equation}
{\mathbb{P}}\big(A_{2^{n}}\big)\leq \exp \bigg(-\frac{\sqrt{\gamma_3}}{\eta \sqrt{n}}%
2^{n}\,\exp \Big(-\frac{\eta ^{2}n}{2 \gamma_3}\Big)\bigg).  \label{Eq:inter-3}
\end{equation}

By taking $\eta >0$ small enough such that $\eta ^{2} <2 \gamma_3 \ln 2$, we
have $\sum_{n=1}^{\infty }{\mathbb{P}}\big(A_{2^{n}}\big)$ $<\infty $. Hence
the Borel-Cantelli lemma implies that almost surely, 
\begin{equation*}
\max_{2\leq k\leq 2^{n}}\big |Z(\x_{n,k})-Z(\x_{n,k-1})\big |\geq \eta
2^{-n\nu /2}\sqrt{n}
\end{equation*}%
for all $n$ large enough. This implies that the right-hand side of (\ref%
{Eq:Umod-lower2}) is bounded from below almost surely by $\eta >0$.
This finishes the proof of Theorem 2. 
 
\bigskip
 
 {\bf Acknowledgements}   The research of Y. Xiao is supported in part by the NSF 
grant DMS-1855185.

\end{document}